\documentclass[dvips,12pt,a4paper]{article}
\usepackage{latexsym,amsmath,amssymb,amscd}
\begin{document}

\thispagestyle{empty}
\baselineskip=6.4mm

\begin{center}
{\bf EPIMORPHISMS, DOMINIONS AND $\mathcal H$-COMMUTATIVE SEMIGROUPS}\\
\end{center}
\begin{center}
{\bf  Noor Alam }\\
\noindent
{\it Department of Mathematics\\
College of Preparatory Year, University of Hail, Hail-2440, K.S.A.\\
email: alam.nooramu@gmail.com}
\end{center}

\begin{center}
{\bf Peter M. Higgins}\\
\noindent
{\it Dept. of Mathematical Sciences\\
University of Essex, Colchester, U.K. CO1 4SQ\\
email: peteh essex.ac.uk}\\
\end{center}

\begin{center}
{\bf Noor Mohammad Khan}\\
\noindent
{\it Department of Mathematics\\
Aligarh Muslim University, Aligarh-202002, India\\
email: nm\_khan123@yahoo.co.in}\\
\end{center}

\noindent
{\bf \large {Abstract}} In the present paper, a series of results and examples that explore the structural features of $\mathcal H$-commutative semigroups are provided. We also generalize a result of Isbell from commutative semigroups to $\mathcal H$-commutative semigroups by showing that the dominion of an $\mathcal H$-commutative semigroup is $\mathcal H$-commutative. We then use this to generalize Howie and Isbell's result that any $\mathcal H$-commutative semigroup satisfying the minimum condition on principal ideals is saturated. \\



\vspace{.3cm}
\noindent{\bf{\large{1. Introduction and Preliminaries}}}

\vspace{.3cm}
In this article,  we are concerned with a series of results and examples that explore the class of semigroups $S$ for which Green's relation $\mathcal H$ is commutative: $ab{\mathcal H}ba~$ for all $a,b$ in $S$. This definition of $\mathcal H$-commutativity was introduced by Tully in \cite {tu}. In (\cite {n}, Theorem 5.1), Nagy proposed a second definition of $\mathcal H$-commutativity: [Def. 1.1]. He then proved that the two characterizations coincide.

\vspace{.3cm}\noindent
{\bf{Definition 1.1}} ([11, Chapter V]) A semigroup $S$ is called $\mathcal H$-commutative if for all $a,b\in S$, there exists $x\in S^{1}$ such that $ab=bxa$.

\vspace{.3cm}
Since such equations are always solvable in any group, we see at once that the collection $HC$ of all $\mathcal H$-commutative semigroups represents an umbrella class for the classes of Groups and Commutative semigroups.

\vspace{.3cm}\noindent
{\bf{Result 1.2}} ([11, Theorem 5.2, Chapter V]) A semigroup $S$ is $\mathcal H$-commutative if and only if Green's equivalence ${\cal H}$ on $S$ is a commutative
congruence on $S$.

\vspace{.3cm}\noindent
{\bf{Result 1.3}} ([11, Theorem 5.3, Chapter V]) Every $\mathcal H$-commutative semigroup is decomposable into a semilattice of archimedean semigroups.

\vspace{.3cm}
We now introduce dominions of semigroups. Dominions of permutative semigroups were studied in \cite {kh1} by Khan and Shah. We restate here the presentation given in \cite {kh1}, adapted to the present context. Let $S$ be any semigroup with a subsemigroup $U$. An element $d\in S$ is said to be \emph {dominated} by $U$ if for every semigroup $T$ and for all homomorphisms $\alpha,\beta:S\to T$, $u\alpha=u\beta$ for all $u \in U$ implies that $d\beta=d\alpha$. The set of all elements of $S$ dominated by $U$ is called the \emph {dominion} of $U$ in $S$ and will be denoted by Dom$(U,S)$. It may be easily checked that Dom$(U,S)$ is a subsemigroup of $S$ containing $U$. Any subsemigroup $U$ of a semigroup $S$ is said to be \emph {closed} in $S$ if Dom$(U,S)$ = $U$ and \emph {absolutely closed} if it is closed in every containing semigroup $S$. Further a semigroup $U$ is said to be \emph {saturated} if Dom$(U,S)\not=S$ for every properly containing semigroup $S$ and \emph {epimorphically embedded or dense} in $S$ if Dom$(U,S)=S$.

\vspace{.3cm}
A (semigroup) morphism $\alpha:S\to T$ is said to be an \emph {epimorphism} (\emph {epi} for short) if for all morphisms $\beta,\gamma$ with domain $T$, $\alpha\beta=\alpha\gamma$ implies $\beta=\gamma$ (where the composition of  morphisms is written from left to right). One may easily check that a morphism $\alpha:S\to T$ is epi if and only if $i:S\alpha\to T$ is epi and the inclusion map $i:U\to S$ is epi if and only if Dom$(U,S)=S$. Every onto morphism is easily seen to be an epimorphism, but the converse is not true in general.

\vspace{.3cm}
Semigroup dominions have been characterized by Isbell's zigzag theorem, which is as follows.

\vspace{.3cm}\noindent
{\bf{Result 1.4}} ([6, Theorem 2.3] or [4, Theorem VIII. 8.3.5]) Let $U$ be a subsemigroup of a semigroup $S$. Then $d\in Dom(U,S)$ if and only if $d\in U$ or there exists a series of factorizations of $d$ as follows:\\
$d~=a_0t_1~=~y_1a_1t_1~=y_1a_2t_2~=~y_2a_3t_2~=\cdots =~y_ma_{2m-1}t_m~=~y_ma_{2m},$
where $m \geq 1,~a_i\in U$ $y_i,t_i\in S\backslash U$  and\\
$~~~~~~~~~~~a_0~=~y_1a_1,~~~~~~~~~~~~a_{2m-1}t_m~=~a_{2m};$\\
$~~~~~~~~a_{2i-1}t_i~=~a_{2i}t_{i+1},~~~~~~~~~~~~~~~~~~~~ y_ia_{2i}~=~y_{i+1}a_{2i+1}$ $(1\leq i\leq m-1).~~~~~~~~~~~~~~~~~~~~~~~~~~~~~~~(1)$\\
Such a series of factorizations is called a \emph {zigzag} in $S$ over $U$ with value $d$, length $m$ and spine $u_0, u_1,....,u_{2m}.$

\vspace{.3cm}
A semigroup $S$ is said to be \emph {permutative} if it satisfies a permutation identity $$x_1x_2\cdots x_n=x_{i_1}x_{i_2}\cdots x_{i_n}~(n \geq 2)$$ for some non-trivial permutation $i$ of the set $\{1,2,\ldots ,n\}$. Further $S$ is called \emph {semicommutative} if $i_1\neq 1$ and $i_n\neq n$ and \emph {left (resp. right) semicommutative} if $i_1\neq 1$ (resp. $i_n\neq n$). Permutative semigroups are not saturated in general because commutative semigroups are not saturated. The infinite monogenic semigroup $\langle a \rangle$ generated by $a$ is not saturated since it is epimorphically embedded in the infinite cyclic group generated by $a$ [4, Chapter VII Exercise 2(i)]. In \cite {ho2}, Howie and Isbell showed that commutative semigroups satisfying the minimum condition on principal ideals were saturated. In \cite {kh} Khan extended this result to semicommutative semigroups and in \cite {kh1}, Khan and Shah extended the theorem to right semicommutative semigroups (see [8, Theorem 2.1]).

\vspace{.3cm}
The class of $\mathcal H$-commutative semigroups had been studied by several authors in one way or the other (see \cite {ma}, \cite {ma1}, \cite {n}, \cite {st} and \cite {tu} for example). The class of ${\mathcal H}$-commutative semigroups was first investigated by Tully \cite {tu}. He studied $\mathcal H$-commutative semigroups with the additional property that each congruence is uniquely determined by its kernel relative to a given element $g$ and claimed that $\mathcal H$-commutative semigroups were precisely those semigroups with Green's relation $\mathcal H$ a commutative congruence.  Nagy \cite {n} described archimedean properties of $\mathcal H$-commutative semigroups and showed that every $\mathcal H$-commutative semigroup is a semilattice of archimedean semigroups. Strecker \cite {st}, then, studied $\mathcal H$-commutative semigroups whose lattice of congruences was a chain. He also proved that a semigroup was an $\mathcal H$-commutative archimedean semigroup with an idempotent if and only if it was the ideal extension of a group by a commutative nilsemigroup. In \cite {ma}, Mary studied semigroups whose set of completely regular elements was an $\mathcal H$-commutative set. Mary also gave equivalent characterizations of the condition that $ab{\mathcal H}ba$ element-wise (for given $a,b\in S$) without assuming the whole semigroup $S$ to be $\mathcal H$-commutative [10, Theorem 2.4] and of $\mathcal H$-commutative regular semigroups [10, Theorem 2.7]. For most of the principal results on $\mathcal H$-commutative semigroups, readers are referred to Nagy's book \cite {n}.

\vspace{.3cm}
Our results naturally fall into three parts, which are presented as Sections 2, 3 and 4 respectively. In Section 2, we show that the structural features of commutative semigroups are mirrored in the class $HC$ in that, for any $S\in HC$, all five Green's relations are equal and correspond to the mutual divisibility
of elements. Moreover $S$ is a semilattice of archimedean components. As a consequence of this, it may be easily deduced that the regular members of $HC$ comprise the class of semilattices of groups. We also show that the $HC$ condition on $S$ is equivalent to the requirement that ${\cal H}$ is a congruence
on $S$ and $S/{\cal H}$ is commutative. Section 3 is devoted to the class $HM$ of $\mathcal H$-commutative monoids and provides some examples and remarks showing the distinction between the classes $HC$ and $HM$.

\vspace{.3cm}
In the last section, we prove that the dominion of a $\mathcal H$-commutative semigroup is $\mathcal H$-commutative; this generalizes Isbell's result, from commutative semigroups to $\mathcal H$-commutative semigroups (see [6, Corollary 2.5]). Finally we show that any $\mathcal H$-commutative semigroup satisfying the minimum condition on principal right ideals is saturated which extends Howie and Isbell's result (see [5, Theorem 3.1]) from commutative semigroups to $\mathcal H$-commutative semigroups.

\vspace{.3cm}
Throughout the remainder of the paper $S$ will denote an $\mathcal H$-commutative semigroup
unless otherwise indicated. Background material and facts on semigroups that are
assumed in what follows can be found in Clifford and Preston \cite {cli}, Higgins \cite {hi} and Howie \cite {ho1} and will be used throughout without explicit mention. For a
comprehensive survey on the topic of semilattice decompositions of
semigroups, there is the text \cite {bo}.\\

\newpage

\noindent{\bf{\large{2. Structure of $\mathcal H$-commutative semigroups}}}

\vspace{.3cm}
First we examine the general character of $\mathcal H$-commutative semigroups.

\vspace{.3cm}
\noindent\textbf{Proposition 2.1}. The idempotents of any $\mathcal H$-commutative
semigroup $S$ are central.

\vspace{.3cm}\noindent
\textbf{Proof}. Let $a\in S$ and $e\in E(S)$. Then, for some $x\in S^{1}$,
we have $ea=axe$, whence $eae=axe^{2}=axe=ea.$ Similarly, for some
$y\in S^{1}$, we have $ae=eya$, whence $eae=e^{2}ya=eya=ae$. Therefore
$ea=eae=ae$. Hence each idempotent commutes with every member of
$S$. ~~~~~~~~~~~~~~~~~~~~~~~~~~~~~~~~~~~~~~~~~~~~~~~~~~~~~~~~~~~~~~~~~~~~~~~~~~~~~~~~~~~~~~~~~~~~~~~~~~~~~~~~~~~$\blacksquare$\\

\noindent\textbf{Proposition 2.2}. Let $S$ be an $\mathcal H$-commutative semigroup. Then $aS=Sa~(\forall a\in S)$.

\vspace{.3cm}\noindent
\textbf{Proof.} Take any $a\in S$. Then $ab\in aS~(b\in S)$. As $S$ is $\mathcal H$-commutative, for some $x\in S^1$, we have $ab=bxa\in Sa$. Thus $aS\subseteq Sa$. By symmetry we also have the reverse inclusion
$Sa\subseteq aS$ whence it follows that $aS=Sa$ for all
$a\in S$, as required. ~~~~~~~~~$\blacksquare$\\

\noindent\textbf{Theorem 2.3}. All five Green's relations coincide on an $\mathcal H$-commutative
semigroup $S$.

\vspace{.3cm}\noindent
\textbf{Proof.} For any $a\in S$, by Proposition 2.2, we have $aS=Sa$ whence it follows that $aS^{1}=S^{1}a$ for all
$a\in S$. Hence, for any $a,b\in S$,
\[
a{\cal L}b\Leftrightarrow(aS^{1}=bS^{1})\Leftrightarrow(S^{1}a=S^{1}b)\Leftrightarrow a{\cal R}b.
\]
So we infer that ${\cal H={\cal L={\cal R}={\cal D}}}$ in $S$. What
is more we have $S^{1}aS^{1}=(S^{1})^{2}a=S^{1}a$, whence it follows
that ${\cal L={\cal J}}$ also and, therefore, all five of Green's
relations coincide on $S$. ~~~~~~~~~~~~~~~~~~~~~~~~~~~~~~~~~~~~~~~~~~~~~~~~~~~~~~~~~~~~~~~~~~~~~~~~~~~~~~~~~~~~~~~~~~~~~~~~~~~~~~~~~~~$\blacksquare$\\

\noindent\textbf{Remark 2.4}. In writing ${\cal H}$, therefore, we have a symbol
that may denote any one of the five Green's relations on $S$, noting
that $a{\cal H}b$ if and only if each of $a$ and $b$ are mutually
\emph{divisible}, meaning that each is a factor of the other. In this
context there is no need to distinguish between left factors, right
factors, or interior factors.

\vspace{.3cm}
Since ${\cal L}$ is a right congruence and ${\cal R}$ is a left
congruence in any semigroup $S$, it follows that in a $\mathcal H$-commutative semigroup,
${\cal H}={\cal L}={\cal R}={\cal D}$ is a congruence.

\vspace{.3cm}
It is the case that if we take any surjective homomorphism $\phi:S\rightarrow T$
from an $\mathcal H$-commutative semigroup $S$, then $T$ is also $\mathcal H$-commutative
since for any $a\phi,b\phi\in T$ $(a,b\in S)$ either $ab=ba$, in
which case $a\phi b\phi=(ab)\phi=(ba)\phi=b\phi a\phi$, or there
exists an $x\in S$ such that $ab=bxa$, in which case $a\phi b\phi=(ab)\phi=(bxa)\phi=b\phi x\phi a\phi$.

\vspace{.3cm}\noindent
We generalize this result in our last section to epimorphisms of $\mathcal H$-commutative semigroups where we show that dominion of any $\mathcal H$-commutative semigroup is $\mathcal H$-commutative i.e. If $U$ is any $\mathcal H$-commutative subsemigroup of a semigroup $S$, then Dom$(U,S)$ is also $\mathcal H$-commutative.

\vspace{.3cm}\noindent
\textbf{Theorem 2.5.}
\begin{enumerate}
\item[(a)] ([11, Theorem 5.2, Chapter V])~ A semigroup $S$ is $\mathcal H$-commutative if
and only if ${\cal H}$ is a congruence and $S/{\cal H}$ is commutative.
\item[(b)] If $S$ is $\mathcal H$-commutative, then $S/{\cal H}$ is the greatest
\emph{combinatorial} (meaning ${\cal H}$-trivial) image of $S$.
\end{enumerate}

\textbf{Proof }(a) By Theorem 2.3 and Remark 2.4, given that $S$ is
$\mathcal H$-commutative, then ${\cal H}$ is a congruence on $S$. Moreover,
for any $a,b\in S$, there exist $x,y\in S^{1}$ such that $ab=bxa=xyba$,
thus showing that $ba|ab$; by symmetry we have $ab|ba$ likewise
so that $ab{\cal H}ba$, whence in $S/{\cal H}$, we have $H_{a}H_{b}=H_{ab}=H_{ba}=H_{b}H_{a}$,
so that $S/{\cal H}$ is commutative.

Conversely suppose that ${\cal H}$ is a congruence on
$S$ and that $S/{\cal H}$ is commutative. Then, for any $a,b\in S$,
we have in $S/{\cal H}$ that $H_{ab}=H_{a}H_{b}=H_{b}H_{a}=H_{ba}$, and $ab{\cal H}ba$.\

\vspace{.3cm}\noindent
(b) Suppose that $H_{a}{\cal H}H_{b}$ in the quotient semigroup $T=S/{\cal H}$.
Then, since $T$ is also $\mathcal H$-commutative, we have $H_{a}|H_{b}$
in $T$ so that, for some $c\in S$, we have $H_{b}=H_{a}H_{c}=H_{ac}$.
Hence $a|ac|b$ so that $a|b$ in $S$ and, by symmetry, $b|a$ also
so that $a{\cal H}b$ in $S$, which is to say that $H_{a}=H_{b}$.
Therefore ${\cal H}$ is trivial in $T$. Hence $T=S/{\cal H}$ is
combinatorial, as required.

\vspace{.3cm}
Conversely let $\mu$ be any congruence on $S$ such that $S/\mu$
is combinatorial. Take any $a,b\in S$ such that $a{\cal H}b$. Then,
since homomorphisms preserve Green's relations, we have $a\mu{\cal H}b\mu$
is $S/\mu$. Since $S/\mu$ is combinatorial, it then follows that
$a\mu=b\mu$, whence we infer that ${\cal H}\subseteq\mu.$ Therefore
${\cal H}$ in the least combinatorial congruence on $S$ (which is
equivalent to saying that $S/{\cal H}$ is the maximum combinatorial
image of $S$), as required.~~~~~~~~~~~~~~~~~~~~~~~~~~~~~~~~~~~~~~~~~~~~~~~~~~~~~~~~~~~~~~~~~~~~~~~~~~~~~~~~~~~~~~~~~~~~~~~~~~~~~~~~~~~~~~~$\blacksquare$

\newpage

\textbf{Theorem 2.6.}
\begin{enumerate}
\item[(a)] For a semigroup $S$, the following are equivalent:
\begin{enumerate}
\item[(i)] $S$ satisfies the equations $\forall a,b\,\exists\,x,y\in S:\,(a=axa)\wedge(ab=bya)$.
\item[(ii)] $S$ is $\mathcal H$-commutative and regular.
\item[(iii)] $S$ is $\mathcal H$-commutative and $S/{\cal H}$ is regular.
\item[(iv)] $S$ is a semilattice of groups.
\item[(v)] ${\cal H=\eta},$ where $\eta$ is the least semilattice congruence
on $S$.
\end{enumerate}
\item[(b)] If $S$ is an $\mathcal H$-commutative semigroup, then,
Reg$(S)$, the set of all regular elements of $S$, if non-empty, is an $\mathcal H$-commutative subsemigroup of $S$ which is
itself a semilattice of groups.
\end{enumerate}
\textbf{Proof} (a) (i) $\Rightarrow$ (ii). The first equation ensures
that $S$ is regular, for then $xax\in V(a)$, while the second ensures
that $S$ is $\mathcal H$-commutative.

\vspace{.3cm}\noindent
(ii) $\Rightarrow$ (iii). By Theorem 2.5, ${\cal H}$ is a congruence
and, so, $S/{\cal H}$ is also regular.

\vspace{.3cm}\noindent
(iii) $\Rightarrow$ (iv). By Theorem 2.5, $S/{\cal H}$ is commutative
and combinatorial and, since $S/{\cal H}$ is also regular, $S/{\cal H}$ consists entirely of idempotents and so $S/{\cal H}$ is a semilattice.
Again, for each $a\in S$, we have $a^{2}\in H_{a}$. Thus $H_{a}$
is a group and $S$ is, therefore, a semilattice of groups.

\vspace{.3cm}\noindent
(iv) $\Rightarrow$ (i). Since $S$ is regular, any $x\in V(a)$ is
a solution to the first equation. Take any $a,b\in S$, whence we
may write $H=H_{ab}=H_{ba}=H_{e}$, where $e$ is the identity element
of the group $\mathcal H$. Hence $be,ea\in H$. Put $y=(be)^{-1}ab(ea)^{-1}$,
where inversion is in the group $H$, noting also that $y\in H$.
Then
\[
bya=b(eye)a=(be)y(ea)=(be)(be)^{-1}ab(ea)^{-1}ea=e(ab)e=ab;
\]
thereby proving that $S$ is an $\mathcal H$-commutative semigroup.

\vspace{.3cm}\noindent
(iv) $\Rightarrow$ (v). In any semigroup, we have ${\cal {\cal J\subseteq\eta}},$
so that ${\cal H\subseteq\eta}$ is always true. (Indeed, since $S$
is regular, we have $\eta={\cal J}^{*}$, the least congruence containing
${\cal J}$). Conversely, since each $\eta$-class $a\eta$ is a group,
it follows that $a\eta\subseteq H_{a}$, a group ${\cal H}$-class,
so that $\eta\subseteq{\cal H}$ and we conclude that ${\cal H=\eta}.$

\vspace{.3cm}\noindent
(v) $\Rightarrow$(iv). Since ${\cal H=\eta}$, it follows that, for
any ${\cal H}$-class $H$ of $S$ and $a\in H$, we have $a^{2}\in H$,
whence $H$ is a group. Since ${\cal H=\eta}$, we have that each
$\eta$-class is a group, and so $S$ is a semilattice of groups.

\vspace{.3cm}\noindent
(b) By Proposition 2.1, idempotents commute with each other whence it follows that Reg$(S)$ is a
subsemigroup as for any $a'\in V(a)$, $b'\in V(b)$ we have $b'a'\in V(ab)$. Again by
Proposition 2.1 it now follows that Reg$(S)$ is a semilattice of groups, whence
from (a) it follows that Reg$(S)$ is an $\mathcal H$-commutative subsemigroup of S.~~~~~~~~~~~~~~~~~~~~~~~~~~~~~~$\blacksquare$

\vspace{.3cm}\noindent
\textbf{Definition 2.7}. A semigroup $S$ is called \emph{archimedean
}if for each $a,b\in S$, there exists $n\geq1$ such that $a|b^{n}$,
both as a left divisor and a right divisor.

\vspace{.3cm}\noindent
\textbf{Remark 2.8}. There is no loss of generality in taking the same
value of $n$ for the left and right divisors, for suppose that $b^{n}=ax$
and $b^{m}=ya$ $(x,y\in S^{1})$. Then $b^{mn}=ax(ax)^{m-1}=(ya)^{n-1}ya$,
so that $a$ is both a left and right divisor of a common power of
$b$.

\vspace{.3cm}\noindent
\textbf{Lemma 2.9}. Let $S$ be an $\mathcal H$-commutative semigroup with
$a,b,a_{1},a_{2},b_{1},b_{2}\in S$ and $n\in\mathbb{Z}^{+}$. Then
\begin{enumerate}
\item[(a)] the relation of divisibility is compatible with multiplication, meaning
that if $a_{1}|b_{1}$ and $a_{2}|b_{2}$, then $a_{1}a_{2}|b_{1}b_{2}$.
\item[(b)] if $a|b$, then $a^{n}|b^{n}$;
\item[(c)] if $a{\cal H}b$, then $a^{n}{\cal H}b^{n}$;
\item[(d)] $(bc)^{n}{\cal H}b^{n}c^{n}$;
\item[(e)] $a|b^{n}c^{n}$ if and only if $a|(bc)^{n}$.
\end{enumerate}
\textbf{Proof. }(a) Since $a_{i}|b_{i}$ $(i=1,2)$, we may write for
some $c_{i}\in S^{1}$ that $b_{i}=c_{i}a_{i}$. Then we have, for
some $x\in S^{1}$, that
\[
b_{1}b_{2}=c_{1}a_{1}c_{2}a_{2}=(c_{1}c_{2}x)(a_{1}a_{2})
\]
 so that $a_{1}a_{2}|b_{1}b_{2}$, as required.

\vspace{.3cm}
(b) Follows by induction on $n$ upon taking $a_{1}=a_{2}=a$ and
$b_{1}=b_{2}=b$.

\vspace{.3cm}
(c) Follows as, by Theorem 2.5, ${\cal H}$ is a congruence on $S$.

\vspace{.3cm}
(d) Follows as, by Theorem 2.5, $S/{\cal H}$ is a commutative semigroup.

\vspace{.3cm}
(e) This follows from (d) and the transitivity of the relation of divisibility.
~~~~~~~~~~~~~~~$\blacksquare$

\vspace{.3cm}
Consider the \emph{archimedean division relation }$\lambda$ on $S$
whereby $a\lambda b$ if $a|b^{n}$ for some $n\ge1$. Clearly $\lambda$
is reflexive. To see that $\lambda$ is transitive, suppose further
that $b|c^{m}$ for some $m\geq1$. Then, by Lemma 2.9(b), we have $a|b^{n}|c^{mn}$,
so that $a\lambda c$, showing that $\lambda$ is transitive. Hence
$\lambda$ is a quasi-order on $S$, which then induces an equivalence
relation $\rho$ on $S$ defined by $a\rho b$ if and only if $a\lambda b$
and $b\lambda a$. Indeed $\rho$ is a congruence on $S$; for suppose that
$a\rho b$ so that $b^{n}=ay$ say, and take any $c\in S^{1}$. Then,
for some $x,y\in S^{1}$, we have $b^{n}c^{n}=ayc^{n}=aycc^{n-1}=acxyc^{n-1}$
so that $ac|b^{n}c^{n}$ whence, by Lemma 2.9(e), we infer that $ac|(bc)^{n}$.
Exchanging the roles of $a$ and $b$ in the argument gives that $bc$
divides some power of $ac$ and so $ac\rho bc.$ Hence $\rho$ is
a right congruence and by the left-right symmetry of the relation
of division, it follows that $\rho$ is also a left congruence and,
therefore, $\rho$ is a congruence on $S$.

\vspace{.3cm}
Recall that for any relation $R\subseteq S\times S$, $R^{*}$ denotes
the least congruence on $S$ that contains the relation $R$. A particular
case of this is that $\eta=\eta_{0}^{*}$, where $\eta_{0}=\{(a,a^{2}),(ab,ba):a,b\in S\}$, is the least semilattice congruence on any semigroup
$S$.
Our discussion has led to the following result, which directly generalizes
the well-known theorem for commutative semigroups {[2, p136]}.

\vspace{.3cm}\noindent
\textbf{Theorem 2.10}. Let $S$ be an $\mathcal H$-commutative semigroup. Then
\begin{enumerate}
\item[(a)] the relation $\rho$ on $S$ defined by $a\rho b$ if and only
if each of $a$ and $b$ divides a power of the other is the least
semilattice congruence $\eta$ on $S$;
\item[(b)] each subsemigroup $a\eta$ $(a\in S)$ of $S$ is archimedean;
\item[(c)] $a\eta$ is a union of ${\cal H}$-classes of $S$ and $a\eta$
contains at most one idempotent.
\end{enumerate}
\textbf{Proof}~(a) From the fact that $a|a^{2}$ and $a^{2}|a^{2}$,
we may conclude that $a\rho a^{2}$ for all $a\in S$. Next $ab=bxa=xyba$
for some $x,y\in S^{1}$ so that $ab|ba$. By symmetry $ba|ab$ and,
so, $ab\rho ba$ for all $a,b\in S$. Since $\eta$ is the least congruence
containing all pairs of the form $(a,a^{2})$ and $(ab,ba)$, it follows
that $\eta\subseteq\rho$ and $\rho$ is itself a semilattice congruence
on $S$.

\vspace{.3cm}
Conversely, suppose that $a\rho b$ so that $a|b^{k}$ and $b|a^{n}$
say. Then, for some $x,y\in S^{1}$, we have $a^{n}=bx$ and $b^{k}=ya$,
which yields:
\[
a\,\eta\,a^{n}=bx\,\eta\,b^{2}x=ba^{n}\,\eta\,ba\,\eta\,b^{k}a=ya^{2}\,\eta\,ya=b^{k}\,\eta\,b
\]
whence it follows that $\rho\subseteq\eta.$ Therefore $\rho=\eta$
as claimed.

\vspace{.3cm}\noindent
(b) Since $S/\eta$ is a band (indeed a semilattice), each class $a\eta$
is a subsemigroup of $S$. Take any $a,b\in S$ such that $a\eta b$
and $x\in S^{1}$ such that $ab=bxa$.

\vspace{.3cm}
We have for some $t\in S^{1}$ and $n\geq1$ that $b^{n}=at$. Then
$b^{2n}=a(tat)\in a\eta$, whence $x=tat\in a\eta$ is such that $ax=b^{2n}$.
By symmetry we conclude that each of $a$ and $b$ divides a power
of the other, on the left and on the right, in the semigroup $a\eta$.
Therefore each $\eta$-class $a\eta$ is an archimedean semigroup.

\vspace{.3cm}\noindent
(c) For any $a\in S$, we have $a{\cal H}b$ if and only if $a|b$
and $b|a$, whence $a\eta b$ so that $H_{a}\subseteq a\eta$. Therefore
each $\eta$-class is a union of Green's classes of $S$. Finally,
if for two idempotents $e,f\in E(S)$, we have $e\,\eta\,f$, then it
follows by idempotency that $e{\cal H}f$ , which implies that $e=f$
(and that $H_{e}$ is the unique maximal subgroup of $S$ contained
in $e\eta$). ~~~~~~~~~~~~~~~~~~$\blacksquare$\\

\vspace{.3cm}
\noindent{\bf{\large{3. More Results and Examples}}}

\vspace{.3cm}
The fact that the $\mathcal H$-commutative condition on a semigroup $S$ is defined
by the first order sentence $(\forall ~a,b\in S, ~\exists ~x\in S), ab=bxa$ ~or~ $ab=ba$ allows us to produce further examples. For the moment,
we first confine ourselves to the class of Monoids.

\vspace{.3cm}\noindent
\textbf{Theorem 3.1}. The class $HM$ of $\mathcal H$-commutative monoids
is closed under the taking of homomorphic images, direct products,
and regular submonoids.

\vspace{.3cm}\noindent
\textbf{Proof}. That $HM$ is closed under the taking of homomorphic
images and direct products follows from $HM$ being defined by the first order sentence
$(\forall~ a,b\in M, ~\exists ~x\in M), ab=bxa$.

\vspace{.3cm}
Next let $N$ be a regular submonoid of $M$. Then, by Theorem 2.6(b), Reg$(M)$
is itself a semilattice of groups, whence the same is true of $N\subseteq$
Reg$(M)$, and, so, by Theorem 2.6(b), $N$ is an $\mathcal H$-commutative
submonoid of $M$. ~~~~~~~~~~~~~~~~~~~~~~~~~~~~~~~~~~~~~~~~~~~~~~~~~~$\blacksquare$

\vspace{.3cm}
The distinction between monoids and semigroups is important. Moreover
$HM$ does not constitute a variety of monoids. Both these conclusions
follow from the following examples. First another closure lemma is proved.

\vspace{.3cm}\noindent
\textbf{Lemma 3.2}. The $0$-direct union $V=S\dot{\cup}T\dot{\cup}\{0\}$
of two $\mathcal H$-commutative semigroups $S$ and $T$\emph{ }is $\mathcal H$-commutative.

\vspace{.3cm}\noindent
\textbf{Proof}. Let $a,b\in V$. If $a,b\in S$ or if $a,b\in T$, then
there exists $x\in S^{1}$ or $x\in T^{1}$ as the case may be such
that $ab=bxa$. Otherwise $ab=ba=0$. So it follows that $V$ satisfies
the condition to be $\mathcal H$-commutative. ~~~~~~~~~~~~~~~~~~~~~~~~~~~~~~~~~~~~~~~~~~~~~~~~~~~~~~~~~~~~~~~~~~~~~~~~~~~~~~~$\blacksquare$

\vspace{.3cm}\noindent
\textbf{Example 3.3}. The class $HC$ is not closed under the taking of direct products
or even under the taking of direct powers of one of its members. What
follows is an example of an $\mathcal H$-commutative semigroup $S$ with
ten elements such that $S\times S$ is not $\mathcal H$-commutative. Let
$C_{3}=\langle a:a^{3}=a^{4}\rangle=\{a,a^{2},a^{3}=z\}$, where $z$
denotes the zero element of this monogenic semigroup. Let $S_{3}$
be the symmetric group on $\{1,2,3\}$ and let $S$ be the $0$-direct
union $S=C_{3}\cup S_{3}\cup\{0\}$. Then $C_{3}$ and $S_{3}$ are
each $\mathcal H$-commutative (as $C_{3}$ is commutative and $S_{3}$
is a group) and, so, by Lemma 3.2, is their $0$-direct union $S$.
We note that $|S|=3+6+1=10$. Take the transpositions $t_{1}=(2$~3)
and $t_{2}=(3\,1)$ of $S_{3}$, noting that $t_{1}t_{2}=(1\,3\,2)\neq(1\,2\,3)=t_{2}t_{1}$.
Consider the product $(a,t_{1})(a,t_{2})$ in $S\times S$ and suppose
that $x\in(S\times S)^{1}$ were such that
$$~~~~~~~~~~~~~~~~~~~~~~~~~~~~~~~~~~~~~~~~~~~~~~~~~~~~~~~(a,t_{1})(a,t_{2})=(a,t_{2})x(a,t_{1}).~~~~~~~~~~~~~~~~~~~~~~~~~~~~~~~~~~~~~~~~~~~~~~~(2)$$

If we had $x=(u,v)\in S\times S$ then, since $aua=z$ for all $u\in C_{3}$,
equation (2) takes the form:
\[
(a,t_{1})(a,t_{2})=(a,t_{2})(u,v)(a,t_{1})
\]
\[
\Rightarrow(a^{2},t_{1}t_{2})=(z,t_{2}vt_{1}),
\]
which is a contradiction as $a^{2}\neq z$. On the other hand if we
put $x=1$, then (2) becomes $(a^{2},t_{1}t_{2})=(a^{2},t_{2}t_{1})$,
which is also false as $t_{1}t_{2}\neq t_{1}t_{1}$. Therefore, although
$S$ is a finite $\mathcal H$-commutative semigroup, $S\times S$ is not
$\mathcal H$-commutative. In particular this shows that the $\mathcal H$-commutative
property cannot be defined by the condition that $S$ satisfies some
set of equations with solutions in $S$ (as opposed to solutions in
$S^{1}$).

\vspace{.3cm}\noindent
\textbf{Example 3.4}. Although $HC$ is closed under the taking of regular subsemigroups,
this is not the case for arbitrary subsemigroups, even if the initial
$\mathcal H$-commutative semigroup happens to be a monoid or indeed a
group. To see this, we need look no further than the free group $G$
on $\{g,h\}$, as $G$ contains the free monoid $M$ on the same
pair of generators and $M$ is clearly not $\mathcal H$-commutative.

\vspace{.3cm}\noindent
\textbf{Example 3.5}. We produce an example of a finite $\mathcal H$-commutative semigroup
$T$ with a subsemigroup $S$ that is not $\mathcal H$-commutative. Let
$U$ be the $0$-disjoint union $U=C_{3}^{1}\cup S_{3}\cup\{0\}$,
so that $U=S\cup\{1\}$, where $S$ is as in Example 3.3. Then $S\times S$
is not $\mathcal H$-commutative and is a subsemigroup of the finite semigroup
$T=U\times U$. What is more, $T$ is $\mathcal H$-commutative: for take
any $(p,q),(u,v)\in T$. We have $(p, q)(u, v) = (pu, qv)$. Let $x = (y, z)$ where, if $p, u \in S_{3}$
we put $y = u^{-1}pup^{-1}$ and otherwise put $y = 1$; similarly put $z = 1$
unless $q, v \in S_{3}$ in which case we put $z = v^{-1}qvq^{-1}$. Then we obtain
$(p, q)(u, v) = (u, v)(y, z)(p, q)$ as $pu = uyp$ and $qv = vzq$, thereby verifying
the $H$-commutative property.

\vspace{.3cm}\noindent
\textbf{Example 3.6}.  The archimedean components of the maximum semilattice image $T/\eta$
of an $\mathcal H$-commutative semigroup $T$ are not necessarily themselves
$\mathcal H$-commutative, as we may verify for the case of the $\mathcal H$-commutative
semigroup $T=U\times U$ of Example 3.5 as follows.

For $i=1,2$, take any $a_{i}=(u_{i},v_{i})\in C_{3}\times S_{3}$.
Then
\[
a_{1}^{3}=(z,v_{1}^{3})=(u_{2},v_{2})(z,v_{2}^{-1}v_{1}^{3})=a_{2}(z,v_{2}^{-1}v_{1}^{3});
\]
from which we conclude that every member $a_{2}\in C_{3}\times S_{3}$
divides some power of every member $a_{1}\in C_{3}\times S_{3}$ in
the semigroup $C_{3}\times S_{3}$.

Next consider any factorization of $a_{1}^{k}$ over $T=U\times U$
of the form $a_{1}^{k}=(a^{k},v_{1}^{k})=(u_{2},v_{2})(u_{3},v_{3})$
$(k\geq1)$. Then $u_{2},u_{3}\in C_{3}^{1}$ and $v_{2},v_{3}\in S_{3}$.
Hence if $a_{2}|a_{1}^{k}$ for some $k\geq1$, then $a_{2}\in C_{3}^{1}\times S_{3}$.
Suppose however that $a_{2}=(1,v_{2})$. Then for any $k\geq1$, we
have $a_{2}^{k}=(1,v_{2}^{k})\not\in a_{1}T$ as $a_{1}\in C_{3}\times S_{3}$.
Therefore we conclude that $C_{3}\times S_{3}$ is an $\eta$-class
of $T$.

As in Example 3.3, we now take $(a,t_{1}),(a,t_{2})\in a_{1}\eta$. However,
since $(a_{1}\eta)^{1}=(C_{3}\times S_{3})^{1}\subseteq(S\times S)^{1}$
and, as shown in Example 3.3, there is no $x\in(S\times S)^{1}$ such that
$(a,t_{1})(a,t_{2})=(a,t_{2})x(a,t_{1})$, it follows that $a_{1}\eta$
is not itself an $\mathcal H$-commutative semigroup.

\vspace{0.8 cm}\noindent{\bf\large{4. Epimorphisms and Dominions}}

\vspace{.3cm}
We now generalize Isbell's result [6, Corollay 2.5]  from commutative semigroups to $\mathcal H$-commutative semigroups.

\vspace{.3cm}\noindent
{\bf Theorem 4.1}. Let $U$ be an $\mathcal H$-commutative subsemigroup of a semigroup $S$. Then Dom$(U,S)$ is $\mathcal H$-commutative.

\vspace{.3cm}\noindent
{\bf Proof}. Let $U$ be any $\mathcal H$-commutative subsemigroup of a semigroup $S$. Then we have to show that Dom$(U,S)$ is also $\mathcal H$-commutative; i.e., for all $d,h\in Dom(U,S)$ there exists some $w\in Dom(U,S)^1$ such that $dh=hwd$.

\vspace{.3cm}\noindent
{\bf Case (i):} If both $d,h\in U$, then, trivially $dh=hwd$ for some $w\in U^1$.

\vspace{.3cm}\noindent
{\bf Case (ii):} Let $d\in U$ and $h\in Dom(U,S)\setminus U$. Then, by the zigzag theorem, there exists a series of factorizations of $h$ as follows:\\
$h~=a_0y_1~=~x_1a_1y_1~=x_1a_2y_2~=~x_2a_3y_2~=\cdots =~x_ma_{2m-1}y_m~=~x_ma_{2m},$
where $m \geq 1,~a_i\in U$ $x_i,y_i\in S\backslash U$  and\\
$~~~~~~~~~~a_0~=~x_1a_1,~~~~~~~~~~~~~~~~~a_{2m-1}y_m~=~a_{2m};$\\
$~~~~~~a_{2i-1}y_i~=~a_{2i}y_{i+1},~~~~~~~~~~~~~~~~~~~~~~~~~~~~~ x_ia_{2i}~=~x_{i+1}a_{2i+1}$ $(1\leq i\leq m-1).~~~~~~~~~~~~~~~~~~~~~~~~(3)$\\

\noindent
Now

\vspace{.1cm}\noindent
$$\begin{array}{lll}
dh&=&da_0y_1~~~~~~~~~~~~~~~~~~~~~~~~~~~~~\mbox{(by zigzag equations (3))}\\
\\
&=&a_0w_1dy_1~~~~~~~~~~~~~~~~~~~~~~~~~~\mbox{(for some $w_1\in U^1$ as $U$ is $\mathcal H$-commutative)}\\
\\
&=&x_1a_1w_1dy_1~~~~~~~~~~~~~~~~~~~~~~~\mbox{(by zigzag equations (3))}\\
\\
&=&x_1w_1dw_2a_1y_1~~~~~~~~~~~~~~~~~~~~\mbox{(for some $w_2\in U^1$ as $U$ is $\mathcal H$-commutative)}\\
\\
&=&x_1w_1dw_2a_2y_2~~~~~~~~~~~~~~~~~~~~\mbox{(by zigzag equations (3))}\\
\\
&=&x_1a_2w_3w_1dw_2y_2~~~~~~~~~~~~~~~~\mbox{(for some $w_3\in U^1$ as $U$ is $\mathcal H$-commutative)}\\
\\
&=&x_2a_3w_3w_1dw_2y_2~~~~~~~~~~~~~~~~\mbox{(by zigzag equations (3))}\\
\\
&=&x_2w_3w_1dw_2w_4a_3y_2~~~~~~~~~~~~~\mbox{(for some $w_4\in U^1$ as $U$ is $\mathcal H$-commutative)}\\
&\vdots&\\
&=&x_mw_{2m-1}w_{2m-3}\cdots w_3w_1 dw_2w_4\cdots w_{2m-2}w_{2m}(a_{2m-1}y_m)\\
\\
&=&x_mw_{2m-1}w_{2m-3}\cdots w_3w_1 dw_2w_4\cdots w_{2m-2}w_{2m}a_{2m}~~~~\mbox{(by zigzag equations (3))}\\
\\
&=&x_ma_{2m}w_{2m+1}w_{2m-1}w_{2m-3}\cdots w_3w_1 dw_2w_4\cdots w_{2m-2}w_{2m}\\&&~~~~~~~~~~~~~~~~~~~~~~~~~~~~~~~~~~~~~~\mbox{(for some $w_{2m+1}\in U^1$ as $U$ is $\mathcal H$-commutative)}
\end{array}$$
$$\begin{array}{lll}
&=&hw_{2m+1}w_{2m-1}w_{2m-3}\cdots w_3w_1w_2w_4\cdots w_{2m-2}w_{2m}w_{2m+2}d\\&&~\mbox{(for some $w_{2m+2}\in U^1$ as $U$ is $\mathcal H$-commutative and by zigzag equations (3))}\\
&=&hwd~\mbox{(where $w = w_{2m+1}w_{2m-1}w_{2m-3}\cdots w_3w_1w_2w_4\cdots w_{2m-2}w_{2m}w_{2m+2}\in U^1$)},
\end{array}$$
as required.\\\\
\noindent
{\bf Case (iii)}: Let $d\in Dom(U,S)\setminus U$ and $h\in U$.\\
\vspace{.5cm}\noindent
The proof in this case is the left-right dual to that of Case (ii).\\
\vspace{.5cm}\noindent
{\bf Case (iv)}: Let $d,h\in Dom(U,S)\setminus U$.\\
\vspace{.5cm}\noindent Let (3) be a zigzag for $h\in Dom(U,S)$ in $S$ over $U$. Now
$$\begin{array}{lll}
dh&=&da_0y_1~~~~~~~~~~~~~~~~~\mbox{(by zigzag equations (3))}\\
\\
&=&a_0w_1dy_1~~~~~~~~~~~~~~\mbox{(by Case(ii) for some $w_1\in U^1$)}\\
\\
&=&x_1a_1w_1dy_1~~~~~~~~~~~\mbox{(by zigzag equations (3))}\\
\\
&=&x_1w_1dw_2a_1y_1~~~~~~~~\mbox{(by Case(ii) for some $w_2\in U^1$ as $w_1d\in Dom(U,S)$)}\\
\\
&=&x_1w_1dw_2a_2y_2~~~~~~~~\mbox{(by zigzag equations (3))}\\
\\
&=&x_1a_2w_3w_1dw_2y_2~~~~~\mbox{(by Case(ii) for some $w_3\in U^1$ as $w_1dw_2\in Dom(U,S)$)}\\
\\
&=&x_2a_3w_3w_1dw_2y_2~~~~~\mbox{(by zigzag equations (3))}\\
\\
&=&x_2w_3w_1dw_2w_4a_3y_2~~\mbox{(by Case(ii) for some $w_4\in U^1$ as $w_3w_1dw_2\in Dom(U,S)$)}\\
&\vdots&\\
&=&x_mw_{2m-1}w_{2m-3}\cdots w_3w_1 dw_2w_4\cdots w_{2m-2}w_{2m}a_{2m-1}y_m\\
\\
&=&x_mw_{2m-1}w_{2m-3}\cdots w_3w_1 dw_2w_4\cdots w_{2m-2}w_{2m}a_{2m}~~~~\mbox{(by zigzag equations (3))}\\
\\
&=&x_ma_{2m}w_{2m+1}w'\\&&{\mbox{(by Case(ii) for some $w_{2m+1}\in U^1$ as $w'\in Dom(U,S)$}},\\~&&~~~~~~~~~~~~~~~~~~~~~~~~~~~ {\mbox {where $w'= w_{2m-1}w_{2m-3}\cdots w_3w_1 dw_2w_4\cdots w_{2m-2}w_{2m}$)}}\\
&=&hw_{2m+1}w_{2m-1}w_{2m-3}\cdots w_3w_1w_2w_4\cdots w_{2m-2}w_{2m}w_{2m+2}d\\&&\mbox{(for some $w_{2m+2}\in U^1$ as $U$ is $\mathcal H$-commutative, }\\&&\mbox {$w'= w_{2m-1}w_{2m-3}\cdots w_3w_1 dw_2w_4\cdots w_{2m-2}w_{2m}$ and by zigzag equations (3))}\\
\\
&=&hwd~\mbox{(where $w = w_{2m+1}w_{2m-1}w_{2m-3}\cdots w_3w_1w_2w_4\cdots w_{2m-2}w_{2m}w_{2m+2}\in U^1$)},
\end{array}$$
as required. Thus Dom$(U,S)$ is $\mathcal H$-commutative. ~~~~~~~~~~~~~~~~~~~~~~~~~~~~~~~~~~~~~~~~~~~~~$\blacksquare$\\

\noindent
{\bf Corollary 4.2}. Let $\phi: S\rightarrow T$ be epi. If $S$ is $\mathcal H$-commutative,
then $T$ is $\mathcal H$-commutative.

\vspace{.3cm}\noindent
{\bf Proof}. As $\phi: S\rightarrow T$ be epi, the inclusion morphism $i:  S\phi\rightarrow T$ is epi. Thus Dom$(S\phi,T)=T$. As $S$ is $\mathcal H$-commutative, by Remark 2.4, $S\phi$ is $\mathcal H$-commutative. Therefore, by Theorem 4.1, $T$ is $\mathcal H$-commutative, as required.~~~~~~~~~~~~~~~~~~~~~~~~~~~~~~~~~~~~~~~~~$\blacksquare$\\

In Propositions 4.3 through 4.7, we assume that $S$ is an $\mathcal H$-commutative semigroup and $T$ is a semigroup containing $S$ such that Dom$(S,T)=T$. We also assume that $K$ is a right ideal of $S$ satisfying the minimum condition on principal right ideals and such that $\forall d\in {T\backslash S},~\forall u\in T$, if $s=du\in S$, then $s\in K$.\\

\noindent
{\bf Proposition 4.3}. For any $a\in K$, there exists $c\in K$ and a positive integer $r$  such that $a^{r}c$ is idempotent and $a^r {\cal H}a^rc$.

\vspace{.3cm}\noindent
{\bf Proof}. Consider the descending sequence $aK^1\supseteq a^2K^1\supseteq...$ of principal right ideals $aK^1, a^2K^1$ etc generated by $a, a^2,...$. By the hypothesis, the above descending sequence must stabilize. Therefore, $$a^rK^1=a^{2r}K^1 ~~\mbox{for~ some}~  r,~ \mbox{and so}~ a^r=a^{2r}c,\eqno{(4)}$$ for some $c\in K$. Then
$$a^r=a^{2r}c=a^ra^rc=a^ra^{2r}cc=a^{3r}c^2=....=a^{(k+1)r}c^k.\eqno{(5)}$$
Hence $a^rc=a^{kr}c^k, ~\mbox {for all}~ k\geq 1$.

\vspace{.3cm}\noindent
Now put $k=2$. We obtain

$$a^rc=a^{2r}c^2=a^{r}a^{r}c^{2}=a^{r}c^{2}xa^r~~ (\mbox {for some}~ x\in S)~~~~~~~~~~~~~~~~~~~~~~~~~~~~~~~~~~~~~~~~~~~~~~~~(6)$$
$$=a^{r}c^{2}xa^{2r}c=a^{r}c^{2}xa^{r}a^{r}c=a^{r}a^{r}c^{2}a^{r}c~~(\mbox {for some}~ x\in S)~~~~~~~~~~~~~~~~~~~~~~~~~~~~~~$$
$$=a^{2r}c^{2}a^{r}c=(a^rc)(a^rc)=(a^rc)^2~~~~~~~~~~~~~~~~~~~~~~~~~~~~~~~~~~~~~~~~~~~~~~~~~~~~~~~~~~~~~$$

\noindent
Therefore, $a^rc$ is an idempotent. Now we show that $a^r{\mathcal {H}}a^rc$. As
$$a^rc=a^{2r}c^2~~~~~~~~~~~~~~~~~~~(\mbox{by equation (6)})~~~~~~~~~~~~~~~~~~~~~~~~~~~~~~~~~~~~~~~~~~~~~~~~~~~~~~~$$
$$=a^{r}a^{r}c^{2}=a^{r}c^{2}xa^r~~~~(\mbox  {for ~some}~x\in S)~~~~~~~~~~~~~~~~~~~~~~~~~~~~~~~~~~~~~~~~~~~~~~~~~~$$
$$=(a^{r}c^{2}x)a^r~~~~~~~~~~~~~~~~~~~~~~~~~~~~~~~~~~~~~~~~~~~~~~~~~~~~~~~~~~~~~~~~~~~~~~~~~~~~~~~~~~~~~~~$$

and $a^r(a^rc)=a^{2r}c=a^r$, we have $a^r{\mathcal {L}}a^rc$. Since, by Theorem 2.3, all Green's relations are equal on $\mathcal H$-commutative semigroups, we have $a^r{\mathcal {H}}a^rc$, as required. ~~~~~~~~~~~~~~~~~$\blacksquare$\\

\noindent
{\bf Proposition 4.4}. For each $b\in T\backslash S$, there exists an idempotent $f\in K$ such that $b=bf~(=fb)$.

\vspace{.3cm}\noindent
{\bf Proof}. As Dom$(S,T)=T$, by the zigzag theorem, $b=a_0x_1=y_1a_1x_1$, for some $a_0, a_1 \in S$. As $y_1a_1=a_0\in S$, by hypothesis, $y_1a_1=a_0\in K$. Hence $b$ has a left divisor $a_0\in K$. Let $B$ be the set of all left divisors of $b$ in $K$. Then $B\not=\emptyset$. Let ${\cal B}$ be the set of all principal right ideals of $K$ generated by the elements of $B$. Let $k\in B$ be such that the principal right ideal of $K$ generated by $k$ is minimal in ${\cal B}$. Then $b=kz$ for some $z\in T\backslash S$. By the same argument used in the factorization of $b$, it follows that $z= k^\prime z^\prime$ for some $z^\prime\in {T\backslash S}$ and $k^\prime\in K$. As the principal right ideal of $K$ generated by $kk^\prime$ is contained in the principal right ideal generated by $k$, we have $k=kk^\prime l=k(k^\prime l)^2=k(k^\prime l)^q$ ($l\in K$) for all $q=1, 2, 3, \ldots$. Now, consider the descending sequence $(k^\prime l)K^1\supseteq (k^\prime l)^2K^1\supseteq...$ of principal right ideals $(k^\prime l)K^1, (k^\prime l)^2K^1$ etc generated by $k^\prime l, (k^\prime l)^2,\ldots$. As $K$ satisfies the minimum condition on principal right ideals, $(k^\prime l)^r=(k^\prime l)^{2r}k^{\prime\prime}$ for some $k^{\prime\prime}\in K$ and some positive integer $r$. Thus, as in the proof of Proposition 4.3, $(k^\prime l)^r$ is a multiple of an idempotent $f=(k^\prime l)^rk^{\prime\prime}$.\\
Hence
$$k=k(k^\prime l)^r=k(k^\prime l)^{2r}k^{\prime\prime}=k(k^\prime l)^r(k^\prime l)^rk^{\prime\prime}=k(k^\prime l)^rf.$$
\noindent
Now
$$b=kz=k(k^\prime l)^rfz=k(k^\prime l)^rfzf~~~~~~~~~~~~~~~~~~\mbox{(by Proposition 2.1)}$$
$$=(kz)f=bf,~~~~~~\mbox{as required.}~~~~~~~~~~~~~~~~~~~~~~~~~~~~~~~~~~~~~~~~~~~~~~~~~~~~~~~~~~~~~~~~~~~~~~~~~{\blacksquare}$$

\vspace{.5cm}
For $b \in T\backslash S$ take $e = f \in K$  as in Proposition 4.4.  Then $b = eb \in eT$.  However $eK \subseteq K \subseteq S$ and $b\not\in S$ so that $eb \in eT\backslash eK$. Hence $eK$  is properly contained in $eT$.\\

\noindent
{\bf Proposition 4.5}. For any idempotent $e\in K$, Dom$(eK,eT)=eT$.

\vspace{.3cm}\noindent
{\bf Proof}. Take any $ed\in eT$ for any $d\in T\backslash S$. Since $d\in Dom(S,T)$, by the zigzag theorem, $d$ has a zigzag in $T$ over $S$. Hence we may write
$$ed=ea_0x_1=(ea_0)(ex_1)~~~~~~~~~~~~~~~~~~~~~~~~~~~(\mbox{by Proposition 2.1})~~~~~~~~~~~~~~~~~~~~~~~~~~~~~~~~~~~~~$$
$$=(ey_1)(ea_1)(ex_1)=(ey_1)(ea_2)(ex_2)~~~~(\mbox{by zigzag equations and Proposition 2.1})$$
$$\begin{array}{lll}
&\vdots~~~~~~~~~~~~~~~~~~~~~~~~~~~~~~~~~~~~~~~~~~~~~~~~~~~~~~~~~~~~~~~~~~~~~~~~~~~~~~~~~~~~~~~~~~~~~~~~~~~~~~~~~~~~\\
\end{array}$$
$$=(ey_m)(ea_{2m})~~~~~~~~~~~~~~~~~~~~~~~~~~~~~~~~~~~~~~~~~~~~~~~~~~~~~~~~~~~~~~~~~~~~~~~~~~~~~~~~~~~~~~~$$

\vspace{.3cm}
which is clearly a zigzag over $eS =eK$ with value $ed$. Therefore, Dom$(eK,eT)=eT$.~~~~~~~~~~~~~~~~~~~~~~~~~~~~~~~~~~~~~~~~~~~~~~~~~~~~~~~~~~~~~~~~~~~~~~~~~~~~~~~~~~~~~~~~~~~~~~~~~~~~~~~~~~~~~$\blacksquare$\\

Recall the natural partial order $\leq$ of the idempotents of a semigroup whereby $e\leq f$~ $\Leftrightarrow$~ $ef=fe=e$. An idempotent $e\in K$ is said to be a right (left) divisor of $y\in {T\backslash S}$ if $y=xe~(y=ex)$ for some $x\in T$. Then necessarily $x\in {T\backslash S}$ as $e\in S$, and $S$ is a subsemigroup of $T$, and since $e$ is an idempotent, we can take $x=y$. By Proposition 2.1, if $S$ is $\mathcal H$-commutative, then $e\leq f$~ $\Leftrightarrow$~$ef=e$ and $e$ is a right divisor of $y$~ $\Leftrightarrow$~$e$ is a left divisor of $y$~ $\Leftrightarrow$~$y=ey=ye$.\\

\noindent
{\bf Proposition 4.6}. For each $b\in T\backslash S$, there exists a smallest idempotent $e\in K$ such that $b=be~(=eb)$.

\vspace{.3cm}\noindent
{\bf Proof}. By Proposition 4.4, there exist $e_1\in E(K)$ such that $b=be_1$. Suppose that $e_2\in E(K)$ is such that $b=be_2$. Then $b=be_2=be_1e_2$ and it
follows that the set $F$ of all idempotent divisors of $b$ in $K$ is a subsemilattice of the semilattice $E(K)$ of all idempotents in $K$. As $K$ satisfies the minimum condition on principal right ideals and idempotents are central, $F$ cannot have  an infinite descending chain and so there must be a least element $e$ in $F$.~~~~~~~~~~~~~~~~~~~~~~~~~~~~~~~~~~~~~~~~~~~~~~~~~~~~~~~~~~~~~~~~~~~~~~~~~~~~~~~~~~~~~~~~~~~~~~~~~~~~~~~~$\blacksquare$\\

Now consider any principal right ideal $B$ of $eK$, for any idempotent $e$, generated by any element $b=ek\in eK$. Then $B=\{b\}\cup b(eK)=\{b\}\cup bK=bK^1$. Therefore $B$ is equal to the principal right ideal of $K$ generated by $b$ whence $eK$, as does $K$, satisfies the minimum condition on principal right ideals. This allows us to apply the argument of Proposition 4.4 to $eK$ in the the following proposition.\\

Moreover $\forall~ d\in eT\backslash eS$, $\forall~ u\in eT$ such that $s=du\in eS$, then we have $d\in T\backslash S$, $u\in T$ and $s=du\in S$. Therefore $s=du\in K$. Now $s=du=edu\in eK$. Thus the conditions that we assumed throughout as regards to factorizations also apply to $eK$ in $eT$ as well.

\vspace{.3cm}\noindent
{\bf Proposition 4.7}. Let $b\in T\backslash S$ and let $e\in E(K)$ be a smallest idempotent such that $b=be$, as provided through Proposition 4.6. If $S$ is properly contained in $T$, then there exist an element $z\in eT\backslash eK$ whose only divisors in $eK$ are the elements of $H_e$, the $\mathcal {H}$-class of the element $e$.

\vspace{.3cm}\noindent
{\bf Proof}. As $b\in T\backslash S$, $b\in eT\backslash eK$. Since Dom$(eK,eT)=eT$, by the zigzag theorem, $b=a_0x_1=y_1a_1x_1$, for some $a_0, a_1\in eK$ and $x_1, y_1\in eT\backslash eK$. As in Proposition 4.4, let $B$ be the set of all divisors of $b$ in $eK$ and let $\cal B$ be the set of all principal right ideals of $eK$ generated by the elements of $B$. Let $k\in B$ be such that the principal right ideal of $eK$ generated by $k$ is minimal in $\cal B$. Then $b=kz$, where $z\in eT\backslash eK$. Let $k'$ be an arbitrary factor of $z$ in $eK$ so that $z = k'z'$  for some $z' \in eT\backslash eK$. Now, as in the proof of Proposition 4.4, there exists $k^* =l({k^\prime{l}})^{n-1}k^{\prime\prime}\in eK$ such that $k^\prime k^*=(k^\prime l)^nk^{\prime\prime}=f$, which is an idempotent in $eK$. This $f$ is an idempotent factor of $z$, and thus also of $b$. Thus $e\leq f$. Hence $k^\prime(k^*e)=fe=e$. As $k^\prime(k^*e)=e$ and $ek^\prime=k^\prime$, we have $k^\prime{\mathcal {R}} e$. Since, by Theorem 2.3, all Green's relations are equal on $\mathcal H$-commutative semigroups, we have $k^\prime {\mathcal {H}}e$. Thus $k^\prime\in H_e$. Hence $z$ is the required element.~~~~~~~~~~~~~~~~~~~~~~~~~~~~~$\blacksquare$

\vspace{.3cm}
\noindent
{\bf Theorem 4.8}. Let $S$ be any $\mathcal H$-commutative semigroup and let $T$ be any semigroup containing $S$ such that Dom$(S,T)=T$. Let $K$ be any right ideal of $S$ satisfying the minimum condition on principal right ideals and such that $\forall d\in {T\backslash S},~\forall u\in T$, if $s=du\in S$, then $s\in K$. Then $S$ is saturated.

\vspace{.3cm}\noindent
{\bf Proof}. Suppose on the contrary that $S$ is not saturated. Then there exists a semigroup $T$ containing $S$ properly and such that Dom$(S,T)=T$. Then, by Propositions 4.5, Dom$(eK,eT)=eT$ for each idempotent $e\in K$. By Proposition 4.7, let $z\in eT\backslash eK$ be such that the only divisors of $z$ in $eK$ are members of $H_e$. Since $T$ (and not just S) is ${\mathcal {H}}$-commutative, the spine factors of $z\in eK$ (for any zigzag of $z\in eT$ over $eK$) are also left(and right) factors of $z$. Therefore the zigzag of $z\in eT$ over $eK$ is in fact a zigzag over $H_e$, a contradiction as $H_e$, being a group, is absolutely closed [5, Theorem 2.3].~~~~~~~~~~~~~~~~~~~~~~~~~~~~~~~~~~~~~~~~~~~~~~~~~~~~~~~~~~~~~~~~~~~~~~~~~~~~~~~~~~~~~~~~~~~~~~~~~~~~~~~~~~$\blacksquare$

\vspace{.5cm}\noindent
Now, if we take $K=S$, then the assumption ``$\forall~ d\in  {T\backslash S},~\forall~ u\in T$, if $s=du\in S$, then $s\in K$" is trivially satisfied. Thus we get the following theorem as a corollary to Theorem 4.8.

\vspace{.3cm}\noindent
{\bf Theorem 4.9}. Any $\mathcal H$-commutative semigroup satisfying the minimum condition on principal right ideals is saturated.\\

We provide a class of examples of a saturated $\mathcal H$-commutative semigroup as an application of Theorem 4.9.

\vspace{.3cm}
\noindent
{\bf Theorem 4.10}. Every $\mathcal H$-commutative archimedean semigroup containing an idempotent element is saturated.

\vspace{.3cm}
\noindent
{\bf Proof}. Let $S$ be any $\mathcal H$-commutative archimedean semigroup containing an idempotent element. By Theorem 4.9, it is sufficient to show that $S$ satisfies the minimal condition on principal right ideals. By [12, Theorem 4], $S$ is an ideal extension of a group $G$ by a commutative nilsemigroup. Now, for any $a\in S$, consider a descending chain $aS^1\supseteq a^2S^1\supseteq \cdots$ of principal right ideals of $S$. Then for some $n\geq 1$, we have
$a^n\in G$ whence $a^{n+k}S^1 = G$ for all $k\geq 0$. In particular the descending
chain stabilizes and the result follows by Theorem 4.9.~~~~~~~~~~~~~~~~~~~~~~~~~~~~~~~~~~~~~~~~~~~~~~~~~~~~~~~~~~~~~~~~~~~~~~~~~~~~~~~$\blacksquare$\\

\vspace{.5cm}\noindent
{{\bf\large Acknowledgement}}.
We sincerely thank the learned referee for his useful and constructive suggestions, including that of Theorem 4.10, that helped considerably to improve the presentation of the paper.\\

\end{document}